\documentclass[11pt,a4paper,sunlimits]{article}
\usepackage{amsmath,amsfonts}
\usepackage{amsthm}
\usepackage{amssymb}
\usepackage{amscd}
\usepackage{xypic}
\usepackage{verbatim}
\usepackage[plainpages=false,colorlinks,hyperindex,pdfpagemode=None,bookmarksopen,linkcolor=red,citecolor=blue,urlcolor=blue]{hyperref}
\usepackage{pdflscape}
\usepackage{stmaryrd}
\usepackage{mathabx}
\usepackage{multirow}
\usepackage{appendix}
 \swapnumbers
\theoremstyle{plain}

\renewcommand{\marginpar}[1] {  }
\renewcommand{\comment}[1] {  }

\usepackage[all]{xy}
\catcode`\é=\active\def é{\'e}
\catcode`\è=\active\def è{\`e}
\catcode`\à=\active\def à{\`a}
\catcode`\ù=\active\def ù{\`u}

\catcode`\ê=\active\def ê{\^{e}}
\catcode`\î=\active\def î{\^{i}}
\catcode`\ô=\active\def ô{\^{o}}
\catcode`\û=\active\def û{\^{u}}
\catcode`\ç=\active\def ç{\c{c}}

\newtheorem{theo}{Theorem}[section]
\newtheorem{lem}[theo]{Lemma}
\newtheorem{prop}[theo]{Proposition}
\newtheorem{cor}[theo]{Corollary}

\theoremstyle{remark}

\numberwithin{equation}{section}

\def \g{\mathfrak}

\def\bb{\backslash}

\def\De{\Delta}
\def\de{\delta}

\def\ga{\gamma}

\def\la{\lambda}
\def\om{\omega}
\def\Om{\Omega}

\def\si{\sigma}

\def\bb{\backslash}

\def\cB{{\mathcal B}}

\def\cE{{\mathcal E}}
\def\cF{{\mathcal F}}

\def\cM{{\mathcal M}}
\def\cO{{\mathcal O}}

\def\cT{{\mathcal T}}

\def\cV{{\mathcal V}}
\def\cW{{\mathcal W}}

\DeclareFontFamily{OT1}{rsfs}{}
\DeclareFontShape{OT1}{rsfs}{n}{it}{<-> rsfs10}{}
\DeclareMathAlphabet{\mathscr}{OT1}{rsfs}{n}{it}

\def\E{{\rm E}}
\def\F{{\rm F}}

\def\C{\mathbb C}

\def\R{\mathbb R}

\def\N{\mathbb N}

\def\R{\mathbb R}

\def\scP{\underline{\mathcal P}}

\def\sA{\underline{ A}}

\def\sG{\underline{ G}}
\def\sH{\underline{ H}}

\def\sJ{\underline{ J}}

\def\sS{\underline{ S}}

\def\me{\medskip}
\def\no{\noindent}
\def\dis{\displaystyle}

\def\ste{\par\smallskip\noindent}
\def\ste{\par\smallskip\noindent}

\def\dem{ {\em Proof~: \ste }}
 \def\beq{\begin{equation}}
\def\eeq{\end{equation}}

\newenvironment{res}
               {\begin{equation}\begin{minipage}{0.85\textwidth}}
               {\end{minipage}\end{equation}}
\def\ber{\begin{res}}
\def\eer{\end{res}}

\def\qed{{\null\hfill\ \raise3pt\hbox{\framebox[0.1in]{}}\break\null}}

\begin{document}

\author{P. Delorme\thanks{The first author was supported by a grant of Agence Nationale de la Recherche with reference
ANR-13-BS01-0012 FERPLAY.}, P. Harinck}
\title{A  local  relative trace formula for $PGL(2)$}
\date{}
\maketitle

\begin{abstract} Following a scheme inspired by B.~Feigon \cite{F}, we describe the  spectral side of a local relative trace formula for $G:= PGL(2, \E)$ relative to the symmetric subgroup $H:=PGL(2,\F)$ where $\E/\F$ is an unramified  quadratic extension of local non archimedean  fields of characteristic $0$. This spectral side is given in terms of regularized normalized periods and normalized $C$-functions of Harish-Chandra. Using the geometric side  obtained in a more general setting by P. Delorme, P. Harinck and S. Souaifi \cite{DHSo}, we deduce a local  relative trace formula for $G$ relative to $H$. We apply our result to invert some orbital integrals.\end{abstract}

 \medskip

\section{Introduction}

Let $\E/\F$ is an unramified  quadratic extension of local non archimedean  fields of characteristic $0$. In this paper,  we prove a local relative trace formula for $G:= PGL(2, \E)$ relative to the symmetric subgroup $H:=PGL(2,\F)$ following a scheme inspired by B.~Feigon \cite{F}.\me

As in \cite{Ar}, the way to establish a local relative trace formula is to describe two asymptotic expansions of a truncated kernel associated to the regular representation of $G\times G$ on $L^2(G)$, the first one in terms of weighted orbital integrals (called the geometric expansion), and the second one  in terms of irreducible representations of $G$ (called the spectral expansion). The truncated kernel we consider is defined as follows. The   regular representation $R$ of $G\times G$ on $L^2(G)$ is given by $(R(g_1,g_2)\psi)(x)=\psi(g_2^{-1} xg_1)$.  For $f=f_1\otimes f_2$, where $f_1$ and $f_2$ are two smooth compactly supported functions on $G$, the corresponding operator $R(f)$ is an integral operator on $L^2(G)$ with smooth kernel
$$K_f(x,y)=\int_G f_1(gy) f_2(xg) dg= \int_G f_1(x^{-1} g y) f_2(g) dg.$$\\
We define the truncated kernel $K^n(f)$ by 
$$K^n(f):=\int_{ H\times H} K_f(x,y) u(x,n) u(y, n) dx dy,$$
where the truncated function $u(\cdot, n)$ is the characteristic function of a large compact subset in $H$ depending on a positive integer $n$ as in \cite{Ar} or \cite{DHSo}. \me

In \cite{DHSo}, we study such a truncated kernel in the more general setting where $H$ is the group of $\F$-points of a reductive algebraic group $\sH$ defined and split over $\F$ and $G$ is the group of $\F$-points of the restriction of scalars $\sG:=\textrm{Res}_{\E/\F}\sH$ from $\E$ to $\F$ and we obtain an asymptotic geometric expansion of this truncated kernel in terms of weighted orbital integrals. \me

It is considerably more difficult to obtain a spectral asymptotic expansion of the truncated kernel and the  main part of this paper is devoted to give it  for $\sH= PGL(2)$. 

First, we express the kernel $K_f$ in terms of  normalized Eisenstein integrals using  the Plancherel formula for $G$ (cf. section \ref{PlanchNorm}) . Then the truncated kernel can be written as   a finite linear combination, depending on unitary irreducible representations  of $G$, of terms involving scalar product of truncated periods (cf. Corollary \ref{SSKn}). The difficulty appears in the terms depending on principal series of $G$. \me

Let $M$ (resp.,  $P$) be the image in $G$ of the group of diagonal (resp., upper triangular) matrices of $GL(2,\E)$ and let  $\bar{P}$ be the parabolic subgroup opposite to $P$. As $M$ is isomorphic to $\E^\times$, we identify characters on $M$ and on $\E^\times$. The group of unramified characters of $M$ is isomorphic to $\C^*$ by a map $z\to\chi_z$. Let  $\de$  be a unitary character of $\E^\times$, which is trivial on a fixed uniformizer of $\F^\times$. For $z\in\C^*$, we set $\de_z:=\de\otimes\chi_z$. We denote by  $(i_P^G\de_z, i_P^G\C_{\de_z})$ the normalized induced representation and by  $(i_P^G\check{\de_z}, i_P^G\check{\C_{\de_z}})$ its contragredient. Then, the normalized truncated period is defined by 
$$P^n_{\de_z}(S) :=\int_H E^0(P,\de_z, S)(h) u(h,n) dh,\quad S\in i_P^G\C_{\de_z}\otimes  i_{\overline{P}}^G\check{\C_{\de_z}},$$
where $E^0(P,\de_z, \cdot)$ is the normalized Eisenstein integral associated to  $i_P^G\de_z$ (cf. (\ref{IE0})).  The contribution of $i_P^G\de_z$ in $K^n(f)$ is a finite linear combination of integrals
$$I_{\de}^n(S,S'):=\int_{\cO} P^n_{\de_z}(S) \overline{P^n_{\de_z}(S') } \frac{dz}{z}, \quad S,S'\in i_P^G\C_{\de_z}\otimes  i_{\overline{P}}^G\check{\C_{\de_z}}$$
where $\cO$ is the torus of complex numbers of modulus equal to $1$. 

To establish the asymptotic expansion of this integral, we recall the notion of normalized regularized period introduced by   B. Feigon (cf. section 4). This period, denoted by 
$$P_{\de_z}(S):=\int_H^* E^0(P,\de_z, S)(h) dh$$
 is meromorphic in a neighborhood $\cV$ of $\cO$ with at most a simple pole at $z=1$ and defines  a $H\times H$ invariant linear form on $ i_P^G\C_{\de_z}\otimes  i_{\overline{P}}^G\check{\C_{\de_z}}$. Moreover, the difference $P_{\de_z}(S)-P^n_{\de_z}(S)$ is a rational function in $z$ on $\cV$ with at most a simple pole at $z=1$ which depends on the normalized $C$-functions of Harish-Chandra. As  normalized Eisenstein integrals and normalized $C$-functions are holomorphic in a neighborhood of $\cO$, we can deduce an asymptotic behavior of the integrals $I_{\de}^n(S,S')$ in terms of   normalized regularized  periods and normalized $C$-functions (cf. Proposition \ref{limites}).\me

Our first result (cf. Theorem \ref{SS}) asserts that $K^n(f)$ is asymptotic to a polynomial function in $n$ of degree $1$ whose coefficients are described in terms of generalized matrix coefficients $m_{\xi,\xi'}$ associated to unitary irreducible representations $(\pi,V_\pi)$ of $G$ where $\xi$ and $\xi'$ are linear forms on  $V_\pi$. When $(\pi, V_\pi)$ is a normalized induced representation, these linear forms are defined from the regularized normalized periods, its residues, and the normalized $C$-functions of Harish-Chandra. \me

We precise the geometric asymptotic expansion of $K^n(f)$ obtained in \cite{DHSo} for $\sH:=PGL(2)$. Therefore, comparing the two asymptotic expansions of $K^n(f)$, we deduce our relative local trace formula and a relation between orbital integrals on elliptic regular points in $H\bb G$ and  some generalized matrix coefficients of induced representations (Theorem \ref{LRTF}).

As corollary of these results, we give an inversion formula for orbital integrals on regular elliptic points of $H\bb G$ and for integral orbitals of a matrix coefficient associated to a cuspidal representation of $G$.

\section{Notation}\label{not}
 Let  $\F$ be a non archimedean local field of characteristic  $0$ and odd  residual characteristic $q$.   Let $\E$ be an unramified quadratic extension of  $ \F$. Let $\cO_\F$ (resp.,  $\cO_\E$) denote the ring of integers in $\F$ (resp., in $\E$). We fix a uniformizer $\om$ in the maximal ideal of $\cO_\F$. Thus $\om$ is also a uniformizer of $\E$. We denote by $v(\cdot)$ the valuation of $\F$, extended to $\E$. Let $|\cdot|_\F$ (resp., $|\cdot|_\E$) denote the normalized valuation   on $\F$ (resp., on $\E$). Thus for $a\in\F^\times$, one has $|a|_\F=|a|_\E^2$.\me

Let $N_{\E/\F}$ be the norm map from $\E^\times$ to $\F^\times$. We denote by $E^1$ the set of elements in $\E^\times$ whose norm is equal to $1$.\me

Let $\sH:= PGL(2)$ defined over $\F$ and let $\sG:=\mbox{Res}_{\E/\F}(\sH\times_\F\E)$ be the restriction of scalars of $\sH$ from $\E$ to $\F$. We set $H:=\sH(\F)=PGL(2,\F)$ and $G:=\sG(\F)= PGL(2,\E)$. Let $K:=\sG(\cO_F)=PGL(2,\cO_E)$.\medskip 

We denote by  $C^\infty(G)$  the space of smooth functions on $G$ and by $C_c^\infty(G)$   the subspace of compactly supported functions in   $C^\infty(G)$. If $V$ is a vector space of valued functions on $G$ which is invariant by right (resp., left) translations, we will denote by $\rho$ (resp., $\la$) the right (resp., left) regular representation of $G$ in $V$.

If $V$ is a vector space, $V'$ will denote its dual. If $V$ is real, 
$V_\C$ will denote its complexification.\me

Let $p$ be the canonical projection of $GL(2,\E)$ onto $G$. We denote by $M$ and $N$ the image   by $p$  of  the subgroups of diagonal matrices and upper triangular unipotent matrices of $GL(2,\E)$ respectively. We set $P:=MN$ and we denote by $\bar{P}$ the parabolic subgroup  opposite to $P$. Let  $\de_P$ be the modular function of $P$. We denote by $1$ and $w$ the representatives in $K$ of the Weyl group $W^G$ of $M$ in $G$. 

For $J=K, M$ or $P$, we set $J_H:= J\cap H$. \me

For $a, b$ in $\E^\times$, we denote by  $diag_G(a,b)$ the image by $p$ of the diagonal matrix $\left(\begin{array}{cc} a&0\\0&b\end{array}\right)\in GL(2,E).$
The natural map $(a,b)\mapsto diag_G(a,b)$ induces an isomorphism from $ E^\times\times E^\times/diag(E^\times)\simeq E^\times$ to $M$ where $diag(E^\times)$ is the diagonal of $ E^\times\times E^\times$. 
\ber\label{car}Hence, each character $\chi$ of $E^\times$ defines a character of $M$ given by $ diag_G(a,b)\mapsto \chi(ab^{-1})$, which we will denote by the same letter. \eer\\
We define the map $h_M:M\to\R$ by   
\beq\label{hM} q^{-h_M(m)}=|ab^{-1}|_\E\quad\textrm{for }  m=diag_G(a,b).\eeq\\
We define similarly $h_{M_H}$ on $M_H$ by $q^{-h_{M_H}(diag_G(a,b))}=|ab^{-1}|_\F$ for $a,b\in F^\times$. Then for $m\in M_H$, one has $\de_P(m)=\de_{P_H}(m)^2=q^{-2h_{M_H}(m)}.$\me

We normalize the Haar measure $dx$ on $F$ so that $\textrm{vol}(\cO_\F)=1$. We define the   measure $d^\times x$ on $\F^\times$ by  $d^\times x=\dfrac{1}{1-q^{-1}}\dfrac{1}{|x|_\F}dx$. Thus, we have $\textrm{vol}(\cO_\F^\times)=1$. 
We let  $M$ and $M_H$ have the measure induced by $d^\times x$. We normalize the Haar measure on $K$ so that $\textrm{vol}(K)=1$. Let $dn$ be the Haar measure on $N$ such  that 
 $$ \int_N\de_{\bar{P}}(m_{\bar{P}}(n)) dn=1.$$ Let $dg$ be the Haar measure on $G$ such that 
 $$\int_G f(g) dg=\int_M\int_N\int_K f(mnk) dk\; dn\; dm.$$\\
We define $dh$ on $H$ similarly. \me

The Cartan decomposition of $H$ is given by \beq\label{Cartan} H=K_H M_H^+K_H \textrm{ where } M_H^+:=\{diag_G(a,b); a, b\in \F^\times, |ab^{-1}|_\F\leq 1\},\eeq
and for any integrable function $f$ on $H$, we have the standard integration formula 
\beq\label{formuleint} \int_H f(x) dx=\int_{K_H}\int_{K_H} \int_{M_H} D_{P_H}(m) f(k_1 m k_2) dm dk_2 dk_1,\eeq
where $$D_{P_H}(m)=\left\{\begin{array}{ll}\de_{P_H}(m)^{-1} (1+q^{-1})& \textrm{ if } m\in M_H^+\\ 0 & \textrm{ otherwise}\end{array}\right.$$
For $h\in H$, we denote by $\cM(h)$ an element of $M_H^+$ such that $h\in K_H\cM(h)K_H$. The element $h_{M_H}(\cM(h))$ is independent of this choice. 
We thank  E. Lapid  who suggests us the proof of the following Lemma.
\begin{lem}\label{hM(mh)} Let $\Om$ be a compact subset of $H$. There is $N_0>0$ sastisfying the following property:\\
 for any $h\in \Om$, there exists $X_h\in\R$ such that, for all $m\in M_H^+$ satisfying $h_{M_H}(m)\geq N_0$, one has $$h_{M_H}(\cM(mh))=h_{M_H}(m)+X_h.$$
\end{lem}
\dem For a matrix $x=(x_{i,j})_{i,j}$ of $GL(2,\F)$, we set
$$F(x):=\log \max_{i,j}\Big( \frac{|x_{i,j}|_\F^2}{|det x|_\F}\Big).$$

The function $F$ is clearly invariant under the action of the center of $GL(2,\F)$, hence it defines a function on $H$ which we denote by the same letter.

Since $|\cdot|_\F$ is ultrametric, for $k\in K_H$ and $h\in H$, we have $F(kh)\leq F(h)$, hence,  $F(k^{-1}kh)\leq F(kh)$. Using the same argument on the right, we deduce that $F$ is right and left invariant by $K_H$.

If $m=diag_G(\om^{n_1}, \om^{n_2})$ with $n_1-n_2\geq 0$ then $F(m)= \log\max\Big(\dfrac{q^{-2n_1}}{q^{-n_1-n_2}}, \dfrac{q^{-2n_2}}{q^{-n_1-n_2}}\Big)=(n_1-n_2)\log q=h_{M_H}(m)\log q.$
Thus, we deduce that  $$F(h)=h_{M_H}(\cM(h))\log q ,\quad h\in H.$$

If $h=p\left(\begin{array}{cc} a & b\\ c& d\end{array}\right)$ and $m=diag_G(\om^{n_1}, \om^{n_2})$, then
$$F(mh)=\log\max\Big(|a|_\F q^{n_2-n_1},|b|_\F q^{n_2-n_1}, |c|_\F q^{n_1-n_2},  |d|_\F q^{n_1-n_2}\Big).$$
Therefore, we can choose $N_0>0$ such that, for any $h\in\Om$ and $m\in M_H^+$ with $h_M(m)>N_0$, we have
$$F(mh)=\log\max\Big(|c|_\F q^{n_1-n_2},|d|_\F q^{n_1-n_2}\Big)=(n_1-n_2)\log q+ \log\max(|c|_\F,|d|_\F).$$
Hence, we obtain the Lemma.\qed
\section{ Normalized Eisenstein integrals and Plancherel formula}\label{PlanchNorm}
We denote by $\widehat{M}_2$ the set of unitary characters of $E^\times$ which are  trivial on $\omega$. For $\de\in \widehat{M}_2$ we let $d(\de)$ be the formal degree of $\de$.

Let $X(M)$ be the complex torus of unramified characters of $M$ and $X(M)_u$ be the compact subtorus of unitary unramified characters of $M$. For $z\in\C^*$, we denote by $\chi_z$ the unramified character of $\E^\times $ defined by $\chi_z(\om)=z$. By definition of $h_M$, we have $\chi_z(m)=z^{h_M(m)/2}$. Each element of $ X(M)$ is of the form $\chi_z$ for some $z\in \C^*$ and $X(M)_u$ identifies with  the group $\cO$ of  complex numbers of modulus  equal to $1$. \\
For  $\de\in \widehat{M}_2$ and  $z\in\C^*$, we set $\de_z:=\de\otimes \chi_z$. We will denote by  $\C_{\de_z}$ the space of $\de_z$. \me

Let $Q=MU$ be equal to $P$ or to $\bar{P}$. Let $\de\in \widehat{M}_2$ and $z\in\C^*$. We denote by $i_Q^G \de_z$ the right representation of $G$ in the space $i_Q^G\C_{\de_z}$ of maps  $v$ from $G$ to $\C$, right invariant by a compact open subgroup of $G$ and such that $v(mug)=\de_Q(m)^{1/2}\de_z(m) f(g)$ for all $m\in M, u\in U$ and $g\in G$.  

One denotes by $(\bar{i}_Q^G\de_z, i_{K\cap Q}^K\C)$ the compact realization of $(i_Q^G \de_z, i_Q^G\C_{\de_z})$ obtained by restriction of functions. If $v\in i_{Q\cap K}^K \C$, one denotes by $v_z$ the element of  $i_Q^G\C_{\de_z}$  whose restriction to $K$ is equal to $v$. \me

One defines a scalar product on $i_{Q\cap K}^K \C$ by 

\beq\label{scalP}(v,v')=\int_K v(k)  \overline{v'(k)} dk,\quad  v,v'\in i_{Q\cap K}^K \C.\eeq

If $z\in\cO$ (hence $\de_z$ is unitary), the representation  $\bar{i}_Q^G(\de_z)$ is unitary. Therefore, by ``transport de structure",  $i_Q^G(\de_z)$ is also unitary.\\
 Let $(\check{\de_z}, \check{\C_{\de_z}})$ be the contragredient representation of $(\de_z,\C_{\de_z})$. We can and will identify  $(i_Q^G\check{\de_z}, i_Q^G\check{\C_{\de_z}})$ with the contragredient representation of $(i_Q^G\de_z,i_Q^G \C_{\de_z})$ and $i_Q^G\C_{\de_z}\otimes i_Q^G\check{\C_{\de_z}}$ with a subspace of $\textrm{End}_G(i_Q^G \C_{\de_z})$ ([W], I.3). \me

Using the isomorphism between   $i_Q^G\C_{\de_z}$ and  $i_{Q\cap K}^K \C$, we can define  the notion of rational or polynomial map from  $X(M)$ to a space depending on $i_Q^G\C_{\de_z}$ as in  (\cite{W} IV.1 and VI.1).\me

 We denote by $A(\bar{Q},Q, \de_z):i_Q^G \C_{\de_z}\to   i_{\bar{Q}}^G\C_{\de_z}$ the standard intertwining operator.  By (\cite{W}, IV. 1. and Proposition IV.2.2),  the map $z\in \C^*\mapsto A(\bar{Q}, Q, \de_z)\in \textrm{Hom}_G(i_Q^G  \C_{\de_z}, i_{\bar{Q}}^G  \C_{\de_z})$ is a rational function on $\C^*$. 
Moreover, there exists a rational complex valued function $j(\de_z)$ depending only on  $M$ such that $A(Q,\bar{Q},\de_z)\circ A(\bar{Q}, Q,\de_z)$ is the dilation of scale $j(\de_z)$. We set
\beq\label{mu} \mu(\de_z):=j(\de_z)^{-1}.\eeq
By (\cite{W} Lemme V.2.1), the map $z\mapsto \mu(\de_z)$ is rational on $\C^*$ and regular on $\cO$.\me

\no The Eistenstein integral $E(Q,\de_z)$ is the map from 
$i_Q^G \C_{\de_z}\otimes i_Q^G\check{\C_{\de_z}}$ to $C^\infty(G)$ defined by
\beq\label{EI}E (Q, \de_z, v\otimes \check{v})(g)=\langle(i_Q^G\de_z)(g)v, \check{v}\rangle,\quad v\in i_Q^G\C_{\de_z}, \check{v}\in i_Q^G\check{\C_{\de_z}}.\eeq \\
 If $\psi\in i_Q^G \C_{\de_z}\otimes i_Q^G\check{\C_{\de_z}}$ is identified with  an endomorphism of $ i_Q^G\C_{\de_z}$,  we have
\beq\label{EItrace}E(P,\de_z,\psi)(g)=\textrm{tr} (i_Q^G\de_z(g)\psi\big).\eeq\\
We introduce the operator   $C_{P,P}(1,\de_z):=Id\otimes A(\bar{P},P,\check{\de_z})$ from $i_P^G \C_{\de_z}\otimes i_P^G\check{\C_{\de_z}}$ to $i_P^G \C_{\de_z}\otimes i_{\bar{P}}^G\check{\C_{\de_z}}$. By (\cite{W}, Lemme V.2.2), one has
\ber\label{adjointC}  the operator  $\mu(\de_z)^{1/2}C_{P,P}(1,\de_z)$ is  unitary and regular on $\cO$.\eer\\
We define the normalized Eisenstein integral $E^0(P,\de_z):  i_P^G \C_{\de_z}\otimes i_{\bar{P}}^G\check{\C_{\de_z}}\to C^\infty(G)$ by
\beq\label{IE0}E^0(P,\de_z, \Psi)=E(P,\de_z, C_{P|P}(1,\de_z)^{-1}\Psi).\eeq
By (\cite{S}, \S 5.3.5) , we have 
\ber\label{IE0holo} $E^0(P,\de_z, \Psi)$ is regular on $\cO$.\eer
 For $f\in C_c^\infty(G)$, we denote by $\check{f}$ the function defined by $\check{f}(g):=f(g^{-1})$. Then, the operator $i_P^G\de_z(\check{f})$ belongs to $ i_P^G \C_{\de_z}\otimes  i_P^G \check{ \C_{\de_z}}\subset \textrm{End}_G( i_P^G \C_{\de_z})$. 
We define the Fourier transform $\cF(P,\de_z, f)\in  i_P^G \C_{\de_z}\otimes i_P^G\check{ \C_{\de_z}}$ of $f$  by
 $$ \cF(P,\de_z, f)=  i_P^G\de_z(\check{f}).$$
 It differs from that of \cite{W} by the constant $d(\de)$.\\
The $G$-invariant scalar product on $i_P^G \C_{\de_z}$ defined in (\ref{scalP}) induces a $G$-invariant scalar product on $i_P^G \C_{\de_z}\otimes i_P^G\check{ \C_{\de_z}}$ given by 
$$(v_1\otimes\check{v_1},v_2\otimes\check{v_2})=(v_1,v_2)(\check{v_1},\check{v_2}).$$
Notice that by the inclusion $i_P^G \C_{\de_z}\otimes i_P^G\check{ \C_{\de_z}}\subset \textrm{End}(i_P^G \C_{\de_z})$, this scalar product coincides with the Hilbert-Schmidt scalar product on the space of Hilbert-Schmidt operators on $ i_P^G\C_{\de_z} $ defined by
\beq\label{scalHS} (S,S')=\textrm{tr}(SS'^*),\eeq\\
where $\textrm{tr}(SS'^*)=\sum_{o.n.b.} \langle SS'^* u_i,u_i\rangle$ and this sum converges absolutely and does not depend on the basis.\\
  Then, the Fourier transform is the unique element of $i_P^G \C_{\de_z}\otimes i_P^G\check{ \C_{\de_z}}$ such that 
\beq\label{adjE}(E(P,\de_z,\Psi), f)_G=(\Psi,  \cF(P,\de_z, f)).\eeq
 Moreover, we have (\cite{W} Lemme VII.1.1)
\beq\label{EFtrace} E(P, \de_z, \cF(P,\de_z,f))(g)=  \textrm{tr} \big[(i_P^G\de_z)(\la(g)\check{f})\big].\eeq\\
We define the normalized Fourier transform $\cF^0(P,\de_z,f)$ of $f\in C_c^\infty(G)$ as the unique element of $i_P^G \C_{\de_z}\otimes i_{\bar{P}}^G\check{\C_{\de_z}}$ such that
$$(\Psi, \cF^0(P,\de_z, f))=(E^0(P,\de_z \Psi), f)_G,\quad\quad \Psi\in  i_P^G \C_{\de_z}\otimes i_{\bar{P}}^G\check{\C_{\de_z}}.$$
It follows easily from  (\ref{adjE}) and (\ref{adjointC}) that 
$$\cF^0(P,\de_z,f)=\mu(\de_z)  C_{P|P}(1,\de_z)\cF(P,\de_z,f),$$
thus we deduce that
 \beq\label{E0E}E^0(P,\de_z,\cF^0(P,\de_z,f))=\mu(\de_z) E(P,  \cF(P,\de_z,f)) .\eeq
Therefore, we can describe the spectral decomposition of the regular representation $R:=\rho\otimes\la$ of $G\times G$ on $L^2(G)$ of  (\cite{W} Théorème VIII.1.1) in terms of normalized Eisenstein integrals as follows. Let $\cE_2(G)$ be the set of  classes of irreducible admissible representations of $G$ whose matrix coefficients are square-integrable.  We will denote by $d(\tau)$ the formal degree of  $\tau\in \cE_2(G)$. Then we have

\beq\label{Plancherel}f(g)=\sum_{\tau\in \cE_2(G)} d(\tau)\textrm{tr}(\tau(\la(g)\check{f}))+\frac{1}{4i\pi}\sum_{\de\in \widehat{M}_2}d(\de) \int_\cO E^0(P,\de_z, \cF^0(P,\de_z,f))(g) \frac{dz}{z}.\eeq

\section{The truncated kernel}\label{Kun}
Let  $f\in C_c^\infty (G\times G)$ be of the form $f(y_1, y_2)= f_1(y_1)f_2(y_2)$ with $f_j\in C_c^\infty(G)$. Then the operator $R(f)$ (where $R:=\rho\otimes \la$) is an integral operator with smooth kernel 
$$K_f(x,y)=\int_G f_1(gy) f_2(xg) dg= \int_G f_1(x^{-1} g y) f_2(g) dg.$$\\
Notice that the kernel studied in \cite{Ar}, \cite{F} or \cite{DHSo} corresponds to the kernel of the representation $\la\times\rho$ which coincides with $K_{f_2\otimes f_1}(x,y)=K_{f_1\otimes f_2}(x^{-1}, y^{-1})$.\me

The aim of this part is to give a spectral expansion of the truncated kernel obtained by integrating $K_f$ against a truncated function on $H\times H$ as in \cite{Ar}.
\begin{lem}\label{Kfspect} For $(\tau, V_\tau)\in\cE_2(G)$, we fix an orthonormal basis $\cB_\tau$ of the space of Hilbert-Schmidt operators on $V_\tau$. For $\de\in\widehat{M}_2$ and $z\in \cO$, we fix  an orthonormal basis $\cB_{\bar{P},P}(\C)$  of  $ i_{P\cap K}^K\C\otimes  i_{\bar{P}\cap K}^K\check{\C}$. Using the isomorphism  $S\mapsto S_z$ between $i_{P\cap K}^K\C\otimes  i_{\bar{P}\cap K}^K\check{\C}$ and $i_{P }^G\C_{\de_z}\otimes  i_{\bar{P} }^G\check{\C_{\de_z}}$, we have 
$$K_f(x,y)=\sum_{\tau\in \cE_2(G)}  \sum_{S\in\cB_{\tau}} d(\tau)\textrm{tr}(\tau(x)\tau(f_1)S\tau(\check{f}_2)) \overline{tr(\tau(y)S)}$$
$$+\frac{1}{4i\pi}\sum_{\de\in \widehat{M}_2}\sum_{S\in\cB_{\bar{P},P}(\C)}d(\de) \int_\cO E^0(P,\de_z, \Pi_{\de_z}(f)S_z)(x)\overline{E^0(P,\de_z,S_z)(y)}\frac{dz}{z},$$
where $\Pi_{\de_z}(f)S_z:=(i_P^G\de_z  \otimes i_{\bar{P}}^G\check{\de_z})(f)S_z=(i_{P}\de_z)(f_1)S_z(i_{\bar{P}}\de_z)(\check{f}_2)$ and the sums over $S$ are all finite.
\end{lem}
\dem For $x\in G$, we set $$h(v):=\int_G f_1(uvx) f_2(xu) du,$$so that \beq\label{Kf}K_f(x,y)=\big[\rho (yx^{-1})h\big](e).\eeq
If $\pi$ is a representation of $G$, one has 
$$\pi\big(\rho (yx^{-1})h\big)=\int_{G\times G} f_1(ugy) f_2(xu)\pi(g) du dg=\int_{G\times G}f_1(u_1)f_2(xu) \pi(u^{-1}u_1y^{-1}) du du_1$$
$$=\int_{G\times G} f_1(u_1) f_2(u_2)  \pi(u_2^{-1}x u_1 y^{-1}) du_1 du_2=\pi(\check{f}_2)\pi(x)\pi(f_1)\pi(y^{-1}).$$ \\
Therefore, using the Hilbert-Schmidt scalar product (\ref{scalHS}), one obtains  for $\tau\in\cE_2(G)$,
\beq\label{trsi} \textrm{tr}\; \tau\big(\rho (yx^{-1})h\big)= \textrm{tr}\;\tau(\check{f}_2)\tau(x)\tau(f_1)\tau(y)^*=(\tau(\check{f}_2)\tau(x)\tau(f_1),\tau(y))$$
$$= \sum_{S\in\cB_\tau}  (\tau(\check{f}_2)\tau(x)\tau(f_1),S^*)\overline{(\tau(y),S^*)} =\sum_{S\in\cB_\tau} \textrm{tr}\;(\tau(x)\tau(f_1)S\tau(\check{f}_2)) \overline{tr(\tau(y)S)},\eeq\\
where the sum over $S$ in $ \cB_{\tau} $ is finite. \me

We consider now  $\pi:= i_P^G\de_z$ with $\de \in\widehat{M}_2$ and $z\in \cO$. 
By (\ref{EFtrace}) and (\ref{E0E}), we have
\beq\label{E0Ftr}E^0(P,\de_z, \cF^0(P,\de_z, \rho (yx^{-1})h))(e)= \mu(\de_z) \textrm{tr } \pi(\rho (yx^{-1})h).\eeq

Let $\cB_{P,P}(\C_{\de_z})$ be an orthonormal basis  of   $ i_{P }^G\C_{\de_z}\otimes  i_{P }^G\check{\C_{\de_z}}$. Since $f_1,f_2\in C_c^\infty(G)$, the operators $\pi(f_1)$ and $\pi(\check{f}_2)$    are of finite rank. Therefore, we deduce as above  that 
$$ \textrm{tr}\; \pi\big(\rho (yx^{-1})h\big)= \textrm{tr}\big( \pi(\check{f}_2)\pi(x)\pi(f_1)\pi(y)^{-1}\big)=\sum_{S\in \cB_{P,P}(\C_{\de_z})}tr (\pi(x)\pi(f_1)S\pi(\check{f}_2))\overline{tr (\pi(y) S)},$$\\
where the sum over $S$ in $ \cB_{P,P}(\C_{\de_z})$ is finite. 

In   what follows, the  sums over elements of an orthonormal basis will be always finite.\me

\no Hence, by (\ref{EItrace}), we deduce that 
\beq\label{tracePi}  \textrm{tr}\;  \pi(\rho (yx^{-1})h)=\sum_{S\in \cB_{P,P}(\C_{\de_z})}E(P,\de_z , \pi(f_1)S\pi(\check{f}_2))(x)\overline{E(P,\de_z,S,y)}. \eeq\\
Recall that we fix  an orthonormal basis $\cB_{\bar{P},P}(\C)$      of  the space $ i_{P\cap K}^K\C\otimes  i_{\bar{P}\cap K}^K\check{\C}$ which is isomorphic to  $ i_{P }^G\C_{\de_z}\otimes  i_{\bar{P} }^G\check{\C_{\de_z}}$ by the map $S\mapsto S_z$.  By (\ref{adjointC}),  the family  $\tilde{S}(\de_z):= \mu(\de_z)^{-1/2} C_{P,P}(1,\de_z)^{-1}S_z$ for $S \in \cB_{\bar{P},P}(\C)$ is an orthonormal basis of   $ i_P^G\C_{\de_z}\otimes i_P^G\check{\C_{\de_z}}$. 

Moreover, using the inclusion $ i_P^G \C_{\de_z}\otimes i_{\bar{P}}^G\check{\C_{\de_z}}\subset \textrm{Hom}_G (i_{\bar{P}}^G\C_{\de_z}, i_P^G\C_{\de_z})$, and the adjonction property of the intertwining operator (\cite{W}, IV.1. (11)),  we have 
$C_{P,P}(1,\de_z)^{-1}S = S\circ A(P,\bar{P},\de_z)^{-1},$
 for all  $S\in  i_P^G \C_{\de_z}\otimes i_{\bar{P}}^G\check{\C_{\de_z}}$.
Since $A(P,\bar{P},\de_z)^{-1}\circ i_{P}^G(\de_z)=i_{\bar{P}}^G(\de_z)\circ A(P,\bar{P},\de_z)^{-1}$, writing  (\ref{tracePi}) for the basis $\tilde{S}(\de_z)$, we obtain
$$ \textrm{tr}\;  \pi(\rho (yx^{-1})h)$$
$$=\mu(\de_z)^{-1 }\sum_{S\in\cB_{\bar{P},P}(\C)}E(P,\de_z,  \pi(f_1)C_{P,P}(1,\de_z)^{-1}(S_z)\pi(\check{f}_2))(x)\overline{E(P,\de_z,C_{P,P}(1,\de_z)^{-1}S_z)}(y)$$
$$=\mu(\de_z)^{-1 }\sum_{S\in\cB_{\bar{P},P}(\C)}E(P,\de_z,  C_{P,P}(1,\de_z)^{-1}[(i_{P}^G\de_z)(f_1)S_z(i_{\bar{P}}^G\de_z)(\check{f}_2)])(x)\overline{E(P,\de_z,C_{P,P}(1,\de_z)^{-1}S_z)}(y)$$
$$=\mu(\de_z)^{-1}\sum_{S\in\cB_{\bar{P},P}(\C)}E^0(P,\de_z,(i_{P}^G\de_z)(f_1)S_z(i_{\bar{P}}^G\de_z)(\check{f}_2))(x)\overline{E^0(P,\de_z,  S_z)(y)}.$$\\
\ber\label{pidef} We  set $\Pi_{\de_z}:=i_P^G\de_z  \otimes i_{\bar{P}}^G\check{\de_z}.$   Then we have   $$\Pi_{\de_z}(f)S_z=(i_{P}^G\de_z)(f_1)S_z(i_{\bar{P}}^G\de_z)(\check{f}_2).$$ \eer\\
By (\ref{E0Ftr}), we obtain $$E^0(P,\de_z,\cF^0(P,\de_z, [\rho(yx^{-1})h]))(e)=\sum_{S\in\cB_{\bar{P},P}(\C)}E^0(P,\pi,\Pi_{\de_z}(f)S_z)(x)\overline{E^0(P,\de_z,S_z)(y)}.$$\\
The Lemma follows from  (\ref{Plancherel}), (\ref{Kf}), (\ref{trsi}) and the above result.\qed


To integrate the kernel $K_f$ on $ H\times H $, we introduce truncation as in \cite{Ar}. Let $n$ be a positive integer. Let  $u(\cdot, n)$ be the truncated function defined on $H$ by
$$u(h,n)=\left\{\begin{array}{cc} 1 & \textrm{ if } h=k_1 m k_2 \textrm{ with }
k_1,k_2\in K_H, m\in H \textrm{ such that } 0\leq |h_{M_H}(m)|\leq n\\
0 & \textrm{otherwise}\end{array}\right.$$
We define the truncated kernel by
\beq\label{Kn} K^n(f):=\int_{ H\times H} K_f(x,y) u(x,n) u(y, n) dx dy.\eeq\\
Since  $K_f(x^{-1},y^{-1})$ coincides with the kernel studied in (\cite{DHSo} 2.2) and $u(x,n)=u(x^{-1},n)$, this definition of the truncated kernel coincides with that of \cite{DHSo}. \\
We  defined truncated periods by 
\beq\label{Pns} P^n_\tau(S) :=\int_H \textrm{tr}(\tau(y)S) u(y,n) dy,\quad   (\tau, V_\tau)\in\cE_2(G), S \in \textrm{End}_{fin. rk}(V_\tau),\eeq
where $\textrm{End}_{fin. rk}(V_\tau)$ is the space of finite rank operators in $\textrm{End}(V_\tau)$, and 
\beq\label{Pnz} P^n_{\de_z}(S) :=\int_H E^0(P,\de_z,S_z)(y)u(y,n) dy, \quad \de\in\widehat{M}_2, z\in\cO, S\in i_{P\cap K}^K\C\otimes i_{\bar{P}\cap K}^K\check{\C}.\eeq

\begin{cor}\label{SSKn} With notation of Lemma  \ref{Kfspect}, one has 
$$K^n(f)=\sum_{\tau\in \cE_2(G)}  \sum_{S\in\cB_{\tau}} d(\tau)P^n_\tau(\tau\otimes\check{\tau}(f)S)\overline{P^n_ \tau(S)}$$
$$+ \frac{1}{4i\pi} \sum_{\de\in \widehat{M}_2}\sum_{S\in\cB_{\bar{P},P}(E)}d(\de) \int_\cO P^n_{\de_z}(\overline{\Pi}_{\de_z}(f)S)\overline{P^n_{\de_z}(S)} \frac{dz}{z},$$
where the sums over $S$ are all finite and $\overline{\Pi}_{\de_z}:=\bar{i}_P^G\de_z\otimes \bar{i}_{\bar{P}}^G\check{\de_z}$.\end{cor}
\dem For $\tau\in \cE_2(G)$ and $S\in\cB_\tau$, one has $\tau(f_1)S\tau(\check{f}_2)=\tau\otimes\check{\tau}(f)S$. Therefore, since the functions we integrate are compactly supported, the assertion follows from Lemma \ref{Kfspect}.\qed
\section{Regularized normalized periods}\label{period}
To determine the asymptotic expansion of the truncated kernel, we recall the notion of regularized period introduced in (\cite{F}). It is defined by meromorphic continuation.

Let $z_0\in\C^*$. Then, for $z\in\C^*$ such that $|zz_0|<1$,  the integral 
$$\int_{M_H^+}   \chi_{z_0}(m)\chi_{z}(m) (1-  u(m,n_0)) dm=\sum_{n>n_0} (zz_0)^n=\frac{(zz_0)^{n_0+1}}{1-zz_0}$$
is well defined and has a meromorphic continuation at $z=1$.  Morever this meromorphic continuation is holomorphic on $\cV-\{ 1\}$ with   a simple pole at $z_0=1$. \me

\no Let $\de\in\widehat{M}_2$. We consider now an holomorphic function  $z\mapsto \varphi_{z }\in C^\infty(G)$ defined in a neighborhood $\cV$ of $\cO$ in $\C^*$ such  that 
\ber\label{phiz} there exist a positive integer $n_0$ and two holomorphic   functions $z\in\cV\mapsto \phi^i_z\in C^\infty(K_H\times K_H), i=1,2$ such that, for $k_1,k_2\in K_H$, and $m\in M_H^+$ satisfying $h_{M_H}(m)> n_0$, we have $$\de_P(m)^{-1/2}\varphi_{z}(k_1mk_2)=\de_{z}(m) \phi^1_z(k_1,k_2)+\de_{z^{-1}}(m) \phi^2_z(k_1,k_2).$$\eer\\
Recall that $\cM(h)$ for $h\in H$ is an element in $M_H^+$ such that $h\in K_H\cM(h) K_H$. By the integral formula (\ref{formuleint}), we deduce that  for $|z|<\min(|z_0|,|z_0|^{-1})$,  the integral 
$$\int_H \varphi_{z_0}(h) \chi_{z}(\cM(h))(1-  u(h,n_0))dh$$
$$=(1+q^{-1})\big(\int_{K_H\times K_H} \phi^1_{z_0}(k_1,k_2)dk_1 dk_2\big)\int_{M_H^+} \de(m) \chi_{z_0z}(m) (1-  u(m,n_0)) dm$$
$$+(1+q^{-1})\big(\int_{K_H\times K_H} \phi^2_{z_0}(k_1,k_2)dk_1 dk_2\big)\int_{M_H^+} \de(m) \chi_{z_0^{-1}z}(m) (1-  u(m,n_0)) dm$$
is also well defined and has a meromorphic continuation at $z=1$.  Morever this meromorphic continuation is holomorphic on $\cV-\{ 1\}$ with  at most a simple pole at $z_0=1$. \me
As $u(\cdot,n_0)$ is compactly supported, we deduce that the integral
$$\int_H \varphi_{z_0}(h) \chi_{z} (\cM(h)) dh= \int_H \varphi_{z_0}(h) \chi_{z}(\cM(h))  u(h,n_0)dh+\int_H \varphi_{z_0}(h) \chi_{z}(\cM(h))(1-  u(h,n_0))dh.$$
 has a meromorphic continuation at $z=1$ which we denote by  
$$\int_H^*\varphi_{z_0}(h) dh.$$
The above discussion implies that $\dis \int_H^*\varphi_{z_0}(h) dh $ is  holomorphic on $\cV-\{1\}$ with at most a simple pole at $z_0=1$.   \me

The next result is established in (\cite{F} Proposition 4.6), but we think that the proof is not complete.  We thank E. Lapid who suggests us the   proof below.
\begin{prop}\label{Hinv}($H$-invariance) For $x\in H$, we have 
$$\int_H^*\varphi_{z_0}(hx) dh=\int_H^*\varphi_{z_0}(h) dh.$$
\end{prop}
\dem We fix $x\in H$. For $z,z'$ in $\C^*$, we set $F(\varphi_{z_0},z,z')(h):= \varphi_{z_0}(h)\chi_z(\cM(h))\chi_{z'}(\cM(hx^{-1}))$. By (\ref{phiz}), for $k_1,k_2\in K_H$, and $m\in M_H^+$ with $h_{M_H}(m)> n_0$, we have
$$ \de_P(m)^{-1/2} F(\varphi_{z_0},z,z')(k_1mk_2)=\phi^1_{z_0}(k_1,k_2)\de(m) (z_0z)^{h_M(m)} z'^{h_M(\cM(k_1mk_2x^{-1}))}$$
$$+\phi^2_{z_0}(k_1,k_2)\de(m) (z_0^{-1}z)^{h_M(m)} z'^{h_M(\cM(k_1mk_2x^{-1}))}.$$
We can choose $n_0$ such that  Lemma \ref{hM(mh)} is satisfied. Thus, for any $k_2\in K_H$, there exists $X_{k_2x^{-1}}\in\R$    such that, for any $m\in M_H^+$ satisfying  $1-u(m,n_0)\neq 0$, we have   $h_{M_H}(\cM(k_1mk_2x^{-1}))=h_{M_H}(m)+X_{k_2x^{-1}}$. We deduce that
$$ \de_P(m)^{-1/2} F(\varphi_{z_0},z,z')(k_1mk_2)(1-u(m,n_0))=\phi^1_{z_0}(k_1,k_2)\de(m) (z_0zz')^{h_{M_H}(m)} z'^{X_{k_2x^{-1}}}$$
$$+\phi^2_{z_0}(k_1,k_2)\de(m) (z_0^{-1}zz')^{h_{M_H}(m)} z'^{X_{k_2x^{-1}}} .$$
Therefore, by Hartogs' Theorem and the same argument as above, the function  $$(z_0,z,z')\mapsto  \int_H \varphi_{z_0}(h)\chi_z(\cM(h))\chi_{z'}(\cM(hx^{-1}))dh $$
is well defined for $|z_0zz'|<1$, and has a meromorphic continuation on $\cV\times (\C^*)^2$.   We denote by $I(\varphi_{z_0},z,z')$ this meromorphic continuation. Moreover,   for $z_0\neq 1$, the function $(z,z')\mapsto I(\varphi_{z_0},z,z')$ is holomorphic in a neighborhood of $(1,1)$.\\
For $|z_0z|<1$, we have $\displaystyle I(\varphi_{z_0}, z, 1)=\int_H \varphi_{z_0}(h)\chi_z(\cM(h)) dh$. Hence we deduce that 
$$ I(\varphi_{z_0}, 1, 1)=\int_H^* \varphi_{z_0}(h) dh.$$
On the other hand, we have  $\displaystyle I(\varphi_{z_0}, 1, z')=\int_H \varphi_{z_0}(hx)\chi_{z'}(\cM(h)) dh$ for $|z_0z'|<1$, therefore, one obtains
$$ I(\varphi_{z_0}, 1, 1)=\int_H^* \varphi_{z_0}(hx) dh.$$
This finishes the proof of the proposition.\qed

We will apply this to   normalized Eisenstein integrals. Let $\de\in\widehat{M}_2$ and $z\in \C^*$. 
Recall that we have defined the operator $C_{P,P}(1,\de_z)$ by $$C_{P,P}(1,\de_z):=Id\otimes A(\bar{P},P,\check{\de_z})\in \textrm{Hom}_G\big( i_{P }^G\C_{\de_z}\otimes  i_{P }^G\check{\C_{\de_z}},  i_{P}^G\C_{\de_z}\otimes  i_{\bar{P}}^G\check{\C_{\de_z}}\big).$$ 
We set  $$C_{P,P}(w,\de_z):=A(P,\bar{P},w\de_z)\la(w)\otimes \la(w)\in \textrm{Hom}_G\big(  i_{P }^G\C_{\de_z}\otimes  i_{P }^G\check{\C_{\de_z}},  i_{P }^G\C_{w\de_z}\otimes  i_{\bar{P} }^G\check{\C_{w\de_z}}\big).$$\\
   where $\la(w)$ is the  left translation by $w$ which induces an isomorphism from $i_{P}^G\C_{\de_z}$ to $i_{\bar{P} }^G\C_{w\de_z} $.    For $ s\in W^G$, we define 
\beq\label{C0}C^0_{P,P}(s,\de_z):=C_{P,P}(s,\de_z)\circ C_{P,P}(1,\de_z)^{-1}\in \textrm{Hom}_G\big(  i_{P }^G\C_{\de_z}\otimes  i_{\bar{P} }^G\check{\C_{\de_z}},  i_{P }^G\C_{s\de_z}\otimes  i_{\bar{P} }^G\check{\C_{s\de_z}}\big).\eeq\\
In particular,   $C^0_{P,P}(1,\de_z)$ is the identity map of $ i_{P }^G\C_{\de_z}\otimes  i_{\bar{P} }^G\check{\C_{\de_z}}$.
By arguments analogous to those of (\cite{W} Lemme V.3.1.), we obtain that 

\ber\label{C0holo} for $s\in W^G$,
  the rational operator $C_{P|P}^0(s,\de_z)$ is regular on $\cO$.\eer \\
  Let   $S\in   i_{P\cap K}^K\C\otimes  i_{\bar{P}\cap K}^K\check{\C}$. By (\ref{IE0holo}), the normalized Eisenstein integral $E^0(P,\de_z,S_z)$ is    holomorphic  in a neighborhood $\cV$ of $\cO$. We may and will assume that $\cV$ is invariant by the map $z\mapsto z^{-1}$. By (\cite{He} Theorem 1.3.1) applied to $\la(k_1^{-1})\rho(k_2)E^0(P,\de_z, S_z), k_1,k_2\in K_H$, there exists a positive integer $n_0$ such that, for $k_1,k_2\in K_H$, and $m\in M_H^+$ satisfying $h_{M_H}(m)> n_0$, we have 
$$\de_P(m)^{-1/2}E^0(P,\de_z, S_z)(k_1mk_2)$$
$$= \de(m)\Big(\chi_z(m) \textrm{tr} \big(\big[C^0_{P,P}(1,\de_z)S_z\big](k_1,k_2)\big)+ \chi_{z^{-1}}(m) \textrm{tr} \big(\big[C^0_{P,P}(w,\de_z)S_z\big](k_1,k_2)\big)\Big).$$
Therefore, the normalized Eisenstein integral  satisfies (\ref{phiz}). Hence, we can define  the normalized regularized period by
\beq\label{regperiod}P_{\de_{z }}(S):=\int_H^*E^0(P,\de_z, S_z)(h) dh, \quad S\in   i_{P\cap K}^K\C\otimes  i_{\bar{P}\cap K}^K\check{\C}.\eeq
The above discussion implies that $P_{\de_{z }}(S)$ is a meromorphic function on the neighborhood $\cV$ of $\cO$ which is  holomorphic on $\cV-\{1\}$.\me

For $s\in W^G$ and $S\in   i_{P\cap K}^K\C\otimes  i_{\bar{P}\cap K}^K\check{\C}$, we set 
\beq\label{fonctionsC}
C(s,\de_z)(S):=(1+q^{-1})\int_{K_H\times K_H}\textrm{tr}\big(\big[C^0_{P,P}(s,\de_z)S_z\big](k_1,k_2)\big)dk_1 dk_2.\eeq

 By the same argument as in (\cite{F} Proposition 4.7), we have the following relations between the truncated period and the normalized regularized period. 
\ber\label{relPnP1} If $\de_{|\F^\times}\neq 1$ then,  for $n$ large enough, we have 
$P_{\de_{z }}(S)=P^n_{\de_{z }}(S),$\eer

\ber\label{relPnP2}  If $\de_{|\F^\times}=1$ then,  for $n$ large enough, we have 
$$P_{\de_{z }}(S)=P^n_{\de_{z }}(S)+  \frac{z^{n+1}}{1-z} C(1,\de_z)(S)+\frac{z^{-(n+1)}}{1-z^{-1}}C(w,\de_z)(S).$$\eer \\
The  following Lemma is   analoguous to  (\cite{F} Lemma 4.8 ). 

\begin{lem}\label{L48} Let $z\in \C^*$ and $S\in   i_{P\cap K}^K\C\otimes  i_{\bar{P}\cap K}^K\check{\C}$. 
\begin{enumerate}\item If $\de_{|F^\times}\neq 1$ and $\de_{|E^1}\neq 1$ then, for $n$ large enough, we have 

$$P_{\de_z}(S)=P^n_{\de_z}(S)=0.$$
\item If $\de_{|F^\times}\neq 1$ and $\de_{|E^1}= 1$ then, for $n$ large enough, we have 
$$P_{\de_z}(S)=P^n_{\de_z}(S).$$
\item If $\de_{|F^\times}= 1$ and $\de_{|E^1}\neq 1$ then $P_{\de_z}(S)=0$ whenever it is defined, and $$C(1,\de_1)(S)= C(w,\de_1)(S).$$
\item  If $\de_{|F^\times}= 1$ and $\de_{|E^1}= 1$ then $\de^2=1$. We have $C(1,\de_1)(S)= -C(w,\de_1)(S)$ and the regularized normalized period $P_{\de_z}(S)$ is meromorphic with  a unique pole  at $z=1$ which is simple.
\end{enumerate}
\end{lem}
\dem Case $2$ follows from (\ref{relPnP1}). By (\cite{JLR} Proposition 22), if $\de_{|E^1}\neq 1$ and $z\neq 1$ then the representation $i_P^G\de_z$ admits no nontrivial $H$-invariant linear form. Thus in that case, Proposition \ref{Hinv} implies $P_{\de_z}(S)=0$ whenever it is defined.  We deduce case $1$ from  (\ref{relPnP1})   and  in case $3$, it follows from (\ref{relPnP2}) that  
$$P^n_{\de_{z }}(S)=-\Big( \frac{z^{n+1}}{1-z} C(1,\de_z)(S)+\frac{z^{-(n+1)}}{1-z^{-1}}C(w,\de_z)(S)\Big).$$
Since $P^n_{\de_{z }}(S)$ and $C(s,\de_z)(S)$ for $s\in W^G$ are  holomorphic functions at $z=1$, and 
\ber\label{ResC}  $$\textrm{Res}( \frac{z^{n+1}}{1-z} C(1,\de_z)(S),z=1) = -  C(1,\de_1)(S),$$
$$\textrm{Res}(\frac{z^{-(n+1)}}{1-z^{-1}}C(w, \de_z)(S),z=1) = C(w,\de_1)(S),$$\eer
we deduce the result in the case 3.\me

\noindent In case 4, we obtain easily  $\de^2=1$. By (\cite{W} Corollaire IV.1.2.), the intertwining operator $A(\bar{P}, P, \de_z)$ has a simple pole at $z=1$. Thus  the function $\mu(\de_z)$ has a zero of order $2$ at $z=1$. In that case, by (\cite{S}, proof of Theorem 5.4.2.1),
  the operators  $C_{P|P}(s,\de_z)$ for $s\in W^G$ have a simple pole at $z=1$ and 
$$\textrm{Res} (C_{P|P}(1,\de_z),z=1)= -\textrm{Res} (C_{P|P}(w, \de_z),z=1).$$
Therefore,  if we set  $T_z:= (z-1) C_{P|P}(1,\de_z)$ and $U_z:= (z-1)C_{P|P}(w, \de_z)$, then  $U_z $ and $T_z^{-1}$ are holomorphic near $z=1$  and $T_1=-U_1$ as $\de^2=1$. By definition (cf. (\ref{C0})), we have $C^0_{P|P}(w,\de_z)= U_zT_z^{-1}$. Hence, one deduces that $C^0_{P|P}(w,\de_1)=-Id=-C^0_{P|P}(1,\de_1)$, where $ Id$ is the identity map of $i_P^G\C_{\de_1}\otimes i_{\bar{P}}^G\check{\C_{\de_1}}$. 
 We deduce the first assertion  in case 4  from the definition of $C(s,\de_z)(S)$ (cf.(\ref{fonctionsC})). 

Since $P^n_{\de_{z }}(S)$ and $C(s,\de_z)(S)$ for $s\in W^G$ are  holomorphic functions at $z=1$, the last assertion follows from   (\ref{relPnP2}), (\ref{ResC}) and the above result. This finishes the proof of the Lemma.\qed

\section{Preliminary  Lemma}
In this part, we prove a preliminary lemma which will allow us to get the asymptotic expansion of the truncated kernel  in terms of regularized normalized periods.\me

Let $\cV$ be a neighborhood of $\cO$ in $\C^*$. We assume that $\cV$ is invariant by the map  $z\mapsto \bar{z}^{-1}$. Let  $f$ be a meromorphic function on $\cV$. We assume that $f$   has  at most a  pole at $z=1$ in $\cV$. \\
For  $r<1$ (resp. $r>1$) such that $f$ is defined on the set of complex numbers of modulus $r$, then the integral $\displaystyle \int_{|z|=r} f(z)dz$ does not depend of the choice of $r$. 
We set
\beq \int_{\cO^-}  f(z)dz:= \int_{|z|=r} f(z)dz,\quad r<1.\eeq
and 
\beq \int_{\cO^+}  f(z)dz:= \int_{|z|=r} f(z)dz,\quad r>1.\eeq\\
Notice that we have \beq\label{residu} \int_{\cO^+} f(z)dz-\int_{\cO^-} f(z)dz=2i\pi \textrm{Res}(f(z),z=1).\eeq\\
The two following properties are easily consequences of the definitions:
\beq\label{limO} \lim_{n\to +\infty}\int_{\cO^-} z^n f(z) dz= 0, \quad\textrm{and} \quad   \lim_{n\to +\infty}\int_{\cO^+} z^{-n} f(z) dz= 0\eeq

We have assumed that $\cV$ is invariant by the map $z\to \bar{z}^{-1}$. Then, the function $\tilde{f}(z):=\overline{f(\bar{z}^{-1})}$ is also a meromorphic function on $\cV$ with at most a pole at $z=1$ and it  satisfies $\tilde{f}(z)=\overline {f(z)}$ for $z\in\cO$.\medskip

Let  $c(s,z)$ and $c'(s,z)$, for $s\in W^G$ be holomorphic functions on $\cV$ such that $c(s,1)\neq 0$ and $c'(s,1)\neq 0$. Let $p$ and $p'$ be two meromorphic functions on $\cV$ with at most a pole at $z=1$. We set 
\ber\label{defpn}$$ p_n(z):= p (z)-\Big[\frac{z^{n+1}}{1-z} c(1,z)+\frac{z^{-(n+1)}}{1-z^{-1}} c(w,z)\Big]$$
and $$p'_n(z):= p' (z)-\Big[\frac{z^{n+1}}{1-z} c'(1,z)+\frac{z^{-(n+1)}}{1-z^{-1}} c'(w,z)\Big].$$\eer
\begin{lem}\label{calcul} We assume that   $p_n$ and $p_n'$ are holomorphic on $\cV$ and that either  $p$ and $p'$ are   vanishing functions or  $c(1,1)=-c(w,1)$ and  $c'(1,1)=-c'(w,1)$ .
Then, the integral 
$$\int_{\cO} p_n(z) \overline{p'_n(z)} \frac{dz}{z}$$
is asymptotic as $n$ approaches $+\infty$ to the sum of 
\beq\int_{\cO^-} \Big(p(z)\tilde{p'}(z) + \frac{c(1,z)\tilde{c'}(1,z)}{(1-z)(1-z^{-1})}+ \frac{c(w,z)\tilde{c'}(w,z)}{(1-z)(1-z^{-1})}\Big) \frac{dz}{z},\eeq
\beq -2i\pi\Big[\frac{d}{dz}\Big(c(w,z)\tilde{c'}(1,z)\Big)\Big]_{z=1}+2i\pi \Big[\frac{d}{dz}\Big(  c(w,z) (z-1)\tilde{p'}(z)+\tilde{c'}(1,z)(z-1)p(z) \Big)\Big]_{z=1},\eeq
and
\beq 2i\pi(2n+1)c(w,1)\tilde{c'}(1,1)-2i\pi (n+1)\big( c(w,1) \textrm{Res}(\tilde{p'},z=1)+\tilde{c'}(1,1) \textrm{Res}(p,z=1) \big).\eeq
\end{lem}
\dem
 Since $ p_n$ and $\tilde{p'}_n$ are holomorphic functions on $\cV$, we have $$\int_{\cO} p_n(z) \overline{p'_n(z)} \frac{dz}{z}=\int_{\cO^-} p_n(z) \tilde{p'_n}(z) \frac{dz}{z}$$
 $$=\int_{\cO^-} \Big( p(z)-\frac{z^{n+1}}{1-z} c(1,z)-\frac{z^{-(n+1)}}{1-z^{-1}} c(w,z)\Big)\Big(\tilde{p'}(z)-\frac{z^{-(n+1)}}{1-z^{-1}} \tilde{c'}(1,z)-\frac{z^{n+1}}{1-z} \tilde{c'}(w,z)\Big) \frac{dz}{z}$$
 $$=\int_{\cO^-} \Big(p(z)\tilde{p'}(z) + \frac{c(1,z)\tilde{c'}(1,z)}{(1-z)(1-z^{-1})}+ \frac{c(w,z)\tilde{c'}(w,z)}{(1-z)(1-z^{-1})}\Big) \frac{dz}{z}$$
 $$+\int_{\cO^-}z^{2(n+1)}\frac{c(1,z)\tilde{c'}(w,z)}{(1-z)^2}\frac{dz}{z}-\int_{\cO^-}z^{n+1}\Big( \frac{c(1,z) \tilde{p'}(z)+p(z) \tilde{c'}(w,z)}{1-z}\Big)\frac{dz}{z}$$
 $$+ \int_{\cO^-}z^{-2(n+1)}\frac{c(w,z)\tilde{c'}(1,z)}{(1-z^{-1})^2}\frac{dz}{z}-\int_{\cO^-}z^{-(n+1)}\Big( \frac{c(w,z) \tilde{p'}(z)+p(z) \tilde{c'}(1,z)}{1-z^{-1}}\Big)\frac{dz}{z}.$$
 
 By (\ref{limO}), the second and third terms of the right hand side converge to $0$ as $n$ approaches   $+\infty$.
 
 By (\ref{residu}), one has 
 $$ \int_{\cO^-}z^{-2(n+1)}\frac{c(w,z)\tilde{c'}(1,z)}{(1-z^{-1})^2}\frac{dz}{z}=\int_{\cO^+}z^{-2(n+1)}\frac{c(w,z)\tilde{c'}(1,z)}{(1-z^{-1})^2}\frac{dz}{z}-2i\pi \textrm{Res}(z^{-2(n+1)}\frac{c(w,z)\tilde{c'}(1,z)}{z(1-z^{-1})^2},z=1).$$
 
 Let $\phi(z):=z^{-2(n+1)}\dfrac{c(w,z)\tilde{c'}(1,z)}{z(1-z^{-1})^2}=z^{-(2n+1)} \dfrac{c(w,z)\tilde{c'}(1,z)}{( z-1)^2}$. Since $c(w,z)$ and $\tilde{c'}(1,z)$ are holomorphic functions on $\cV$, the function  $\phi$ has a unique pole of order $2$ at $z=1$. Thus, we obtain
 $$\textrm{Res}(\phi,z=1)= \Big[\frac{d}{dz}\Big( (z-1)^2 \phi(z)\Big)\Big]_{z=1}=-(2n+1)c(w,1)\tilde{c'}(1,1)+\Big[\frac{d}{dz}\Big(c(w,z)\tilde{c'}(1,z)\Big)\Big]_{z=1}.$$
 We deduce from (\ref{limO}) that 
\beq\label{I1} \int_{\cO^-}z^{-2(n+1)}\frac{c(w,z)\tilde{c'}(1,z)}{(1-z^{-1})^2}\frac{dz}{z}=2i\pi(2n+1)c(w,1)\tilde{c'}(1,1)-2i\pi\Big[\frac{d}{dz}\Big(c(w,z)\tilde{c'}(1,z)\Big)\Big]_{z=1}+\epsilon_1(n),\eeq
 where $\displaystyle\lim_{n\to+\infty}\epsilon_1(n)=0$.\me

When $p$ and $p'$ are vanishing functions, we obtain the result of the Lemma.\\
Otherwise, by (\ref{defpn}) and  our assumptions, the function $\dfrac{c(w,z) \tilde{p'}(z)+p(z) \tilde{c'}(1,z)}{1-z^{-1}}$ is  a meromorphic function with a unique pole of order $2$ at $z=1$. Applying  the same argument as above, we obtain
 $$\int_{\cO^-}z^{-(n+1)}\Big( \frac{c(w,z) \tilde{p'}(z)+p(z) \tilde{c'}(1,z)}{1-z^{-1}}\Big)\frac{dz}{z}$$
 $$=\int_{\cO^+}z^{-(n+1)}\Big( \frac{c(w,z) \tilde{p'}(z)+p(z) \tilde{c'}(1,z)}{1-z^{-1}}\Big)\frac{dz}{z}-2i\pi \Big[\frac{d}{dz}\Big( z^{-(n+1)}(z-1)(  c(w,z) \tilde{p'}(z)+p(z) \tilde{c'}(1,z))\Big)\Big]_{z=1}$$
$$=2i\pi (n+1)\big( c(w,1) \textrm{Res}(\tilde{p'},z=1) + \textrm{Res}(p ,z=1)\tilde{c'}(1,1)\big)$$ 
$$-2i\pi \Big[\frac{d}{dz}\Big(  c(w,z) (z-1)\tilde{p'}(z)+(z-1)p(z) \tilde{c'}(1,z)\Big)\Big]_{z=1}+\epsilon_2(n),$$
 where $\displaystyle\lim_{n\to+\infty}\epsilon_2(n)=0$. \me

 Therefore, we obtain the Lemma by (\ref{I1}) and the above result.\qed
\section{Spectral side of a local relative trace formula}

We recall the spectral expression of the truncated kernel obtained in Corollary  \ref{SSKn}:
 $$K^n(f)=\sum_{\tau\in \cE_2(G)}  \sum_{S\in\cB_{\tau}} d(\tau)P^n_ \tau(\tau\otimes\check{\tau}(f)S)\overline{P^n_\tau(S)}$$
$$+ \frac{1}{4i\pi} \sum_{\de\in \widehat{M}_2}\sum_{S\in\cB_{\bar{P},P}(E)}d(\de) \int_\cO P^n_{\de_z}(\overline{\Pi}_{\de_z}(f)S)\overline{P^n_{\de_z}(S)} \frac{dz}{z},$$
where the sums over $S$ are all finite and $\overline{\Pi}_{\de_z}:=\bar{i}_P^G\de_z\otimes \bar{i}_{\bar{P}}^G\check{\de_z}$.\\
By (\cite{F} Lemma 4.10), if $(\tau, V_\tau)\in \cE_2(G)$ and $S\in\textrm{End}_{fin.rk}(V_\tau)$, then
\beq\label{limPns}\lim_{n\to +\infty} P_\tau^n(S)=\int_H \textrm{tr}(\tau(h)S) dh.\eeq
We   consider now  the second term of  the above expression of $K^n(f)$. 
Let $\de\in\widehat{M}_2$ and $S\in   i_{P\cap K}^K\C\otimes  i_{\bar{P}\cap K}^K\check{\C}$. We keep notation of the previous section. In particular, for $z\in \C^*$, we have $\tilde{C}(s,\de_z)(S)=\overline{C(s,\de_{\bar{z}^{-1}})(S)}$ and $\tilde{P}_{\de_z}(S)=\overline{P_{\bar{z}^{-1}}(S)}.$  By definition of $\de_z$, we have  $\de_1=\de$.
\begin{prop}\label{limites} Let $S\in   i_{P\cap K}^K\C\otimes  i_{\bar{P}\cap K}^K\check{\C}$. We set $S'_z:=\overline{\Pi}_{\de_z}(f)S$.
\begin{enumerate}
\item  If $\de_{|\F^\times}\neq 1$ and $\de_{|\E^1}\neq 1$ then, for $n\in\N$ large enough, one has $$\int_\cO P^n_{\de_z}(S'_z)\overline{P^n_{\de_z}(S )} \frac{dz}{z}=0.$$ 
\item  If $\de_{|F^\times}\neq 1$ and $\de_{|E^1}= 1$ then 
$$\lim_{n\to+\infty} \int_\cO P^n_{\de_z}(S'_z)\overline{P^n_{\de_z}(S)}\frac{dz}{z}= \int_\cO P_{\de_z}(S'_z)\overline{P_{\de_z}(S )}\frac{dz}{z}.$$ 

\item  Assume that  $\de_{|F^\times}= 1$ and $\de_{|E^1}\neq 1$.  Then   $$ \int_\cO P^n_{\de_z}(S'_z)\overline{P^n_{\de_z}(S )} \frac{dz}{z}$$
is asymptotic when $n$ approaches $+\infty$ to $$2i\pi(2n+1)C(1,\de)(S'_1)\overline{C(1,\de)(S)}$$
$$+\int_{\cO^-} \Big(  \frac{C(1,\de_z)(S'_z)\tilde{C}(1,\de_z)(S)}{(1-z)(1-z^{-1})}+ \frac{C(w,\de_z)(S'_z)\tilde{C}(w,\de_z)(S)}{(1-z)(1-z^{-1})}\Big) \frac{dz}{z}$$
$$-2i\pi\frac{d}{dz}\Big[C(w,\de_z)(S'_z) \tilde{C}(1,\de_z)(S) \Big]_{z=1}.$$
\item  Assume that $\de_{|F^\times}= 1$ and $\de_{|E^1}= 1$.  Then
$$ \int_\cO P^n_{\de_z}(S'_z)\overline{P^n_{\de_z}(S)} \frac{dz}{z}$$
is asymptotic when $n$ approaches   $+\infty$ to
$$ 2i\pi(2n+3)C(1,\de)(S'_1)\overline{C(1,\de)(S)}$$
$$+\int_{\cO^-} \Big( P_{\de_z}(S'_z)\overline{P_{\de_z}(S)} + \frac{C(1,\de_z)(S'_z)\tilde{C}(1,\de_z)(S)}{(1-z)(1-z^{-1})}+ \frac{C(w,\de_z)(S'_z)\tilde{C}(w,\de_z)(S)}{(1-z)(1-z^{-1})}\Big) \frac{dz}{z}$$
$$-2i\pi\frac{d}{dz}\Big[C(w,\de_z)(S'_z) \tilde{C}(1,\de_z)(S) \Big]_{z=1}$$
$$+2i\pi \Big[\frac{d}{dz}\Big((z-1)P_{\de_z}(S'_z) \tilde{C}(1,\de_z)(S)+C(w,\de_z)(S'_z) (z-1)\tilde{P}_{\de_z}(S)\Big) \Big]_{z=1}.$$
\end{enumerate}

\end{prop}
\proof The two first assertions  are immediate consequences of Lemma \ref{L48}. To prove {\it 3.} and {4.}, we set:\\
$$p_n(z):=P^n_{\de_z}(S'_z(f)), \quad p'_n(z):= P^n_{\de_z}(S),\quad\quad p(z):=  P_{\de_z}(S'_z(f)), \quad p'(z):= P_{\de_z}(S)$$\\
and $c(s,z):=   C(s,\de_z)(S'_z(f)), \quad c'(s,z):=   C(s,\de_z)(S)$ for $s\in W^G .$\me

\no By (\ref{relPnP2}) and Lemma \ref{L48}, these functions satisfy (\ref{defpn}) and we can apply Lemma \ref{calcul}. The result in case 
  {\it 3}  follows immediately since   $p(z)=p'(z)=0$ by   Lemma \ref{L48}.\me

\no In case {\it 4},  we have $c(1,1)=-c(w,1)$ and $c'(1,1)=-c'(w,1)$ by Lemma \ref{L48}. Moreover, the relations (\ref{defpn}) give
$\textrm{Res}(p,z=1)=-c(1,1)+c(w,1)$ and $\textrm{Res}(\tilde{p'}, z=1)=c'(1,1)-c'(w,1)$. Hence, we obtain
$$2i\pi(2n+1)c(w,1)\tilde{c'}(1,1)-2i\pi (n+1)\big( c(w,1) \textrm{Res}(\tilde{p'},z=1)+\tilde{c'}(1,1) \textrm{Res}(p,z=1) \big).$$
$$=2i\pi(2n+3)c(1,1)\tilde{c'}(1,1),$$
and the result in that case follows from Lemma \ref{calcul}.\qed

To describe the spectral side of our local relative trace formula, we introduce generalized matrix coefficients.

\no Let $(\pi,V)$ be a smooth unitary representation of $G$. We denote by $(\pi',V')$ its  dual representation. Let $\xi$ and $\xi'$ be two linear forms on $V$. For $f\in C_c^\infty(G)$, the linear form $\pi'(\check{f})\xi$ belongs to the smooth dual $\check{V}$ of $V$ (\cite{R} Théorème III.3.4 and I.1.2). The scalar product on $V$ induces an isomorphism $j:v\mapsto (\cdot,v)$ from the conjugate complex vector space $\overline{V}$ of $V$ and $\check{V}$, which intertwines the complex conjugate of $\pi$  and $\check{\pi}$ as $\pi$ is unitary. One has
$$\check{v}(v)=(v, j^{-1}(\check{v})),\quad v\in V, \check{v}\in \check{V}.$$
Therefore, for $v\in V$, we have
$$\big(\pi'(\check{f})\xi\big)(v)=\xi\big(\pi(f)v\big)=(v,j^{-1}\big(\pi'(\check{f})\xi\big)).$$
As $\pi(f)$ is an operator of finite rank, we have for any orthonormal basis $\cB$ of $V$
\beq\label{j-}j^{-1}\big(\pi'(\check{f})\xi\big)=\sum_{v\in\cB}(\pi'(\check{f})\xi)(v)\cdot v\eeq
where the sum over $v$ is finite, and $(\la,v)\mapsto\la\cdot v$ is the action of $\C$ on $\overline{V}$.\\
 Let $\overline{\xi'}$ be   the linear form on $\overline{V}$ defined by
$\overline{\xi'}(u)=\overline{\xi'(u)}.$ We define the generalized matrix coefficient $m_{\xi,\xi'}$ by
$$m_{\xi,\xi'}(f)= \overline{\xi'}\Big(j^{-1}\big(\pi'(\check{f})\xi\big)\Big).$$
Then, by (\ref{j-}), we obtain
\beq\label{matcoef}m_{\xi,\xi'}(f)=\sum_{v\in\cB} \xi(\pi(f)v)\overline{\xi'}(v).\eeq
Hence, this sum is independent of the choice of the basis $\cB$.

Let $z\in \C^*$. We set $(\Pi_z,V_z):=( i_P^G\de_z \otimes i_{\bar{P}}^G\Check{\de_{\de_z}}, i_P^G\C_{\de_z}\otimes i_{\bar{P}}^G\Check{\C_{\de_z}})$. We denote by $(\overline{\Pi}_z, V)$ its compact realization. We define meromorphic linear forms on $V_z$ using the isomorphism $V_z\simeq V$. 
\begin{lem} Let $\xi_z$ and $\xi'_z$ be two linear forms on $V$ which are  meromorphic in $z$  on a neighborhood $\cV$ of $\cO$. Let $\cB$ be an orthonormal basis of $V$. Then, for $f\in C_c^\infty(G\times G)$, the sum
$$\sum_{S\in\cB} \xi_z(\overline{\Pi}_z(f)S) \overline{\xi_{\bar{z}^{-1}}(S)}$$
is a finite sum over $S$ which is independent  of the choice of the basis $\cB$.
\end{lem}

\dem For $z\in\cO$, the representation $\Pi_z$ is unitary. Hence (\ref{matcoef}) gives the Lemma in that case.  
 Since the linear forms $\xi_z$ and $\xi'_z$ are meromorphic on $\cV$, we deduce the result  of the Lemma for any $z$ in $\cV$ by meromorphic continuation.\qed

With notation of the Lemma, we define,  for 
$z\in \cV$, the generalized matrix coefficient $m_{\xi_z,\xi'_{\bar{z}^{-1}}}$ associated to $(\xi_z,\xi'_z)$ by
$$m_{\xi_z,\xi'_{\bar{z}^{-1}}}(f):=\sum_{S\in\cB} \xi_z(\overline{\Pi}_z(f)S) \overline{\xi_{\bar{z}^{-1}}(S)}.$$\\
Therefore, using Proposition \ref{limites}, we can deduce the asymptotic behavior of the truncated kernel in terms of generalized matrix coefficients.
\begin{theo}\label{SS} As $n$ approaches $+\infty$,  the truncated kernel $K^n(f)$ is asymptotic to 
$$n\sum_{\de\in \widehat{M}_2,  \de_{|\F^\times}= 1}d(\de)m_{C(1,\de),C(1,\de)}(f)$$
$$+\sum_{\tau\in \cE_2(G)} d(\tau) m_{P_\tau,P_\tau}(f)+ \frac{1}{4i\pi} \sum_{ \de\in \widehat{M}_2, \de_{|\F^\times}\neq 1, \de_{|\E^1}=1}d(\de) \int_\cO  m_{P_{\de_z},P_{\de_z}}(f)\frac{dz}{z}$$
$$+\frac{1}{4i\pi} \sum_{\de\in \widehat{M}_2, \de_{|\F^\times}= 1} R_{\de}(f)+d(\de)\int_{\cO^-} \frac{ m_{C(1,\de_z), C(1,\de_{\bar{z}^{-1}})}(f)+ m_{C(w,\de_z), C(w,\de_{\bar{z}^{-1}})}(f)}{(1-z)(1-z^{-1})}\frac{dz}{z}$$
$$+\frac{1}{4i\pi} \sum_{\de\in \widehat{M}_2, \de_{|\F^\times}=\de_{|E^1}=1}R_\de(f)+\tilde{R}_{\de}(f)+d(\de)\int_{\cO^-}m_{P_{\de_z},P_{\de_{\bar{z}^{-1}}}}(f)\frac{dz}{z}.$$
where
$$R_{\de}(f):=2i\pi d(\de)\Big(m_{C(1,\de),C(1,\de)}(f)- \Big[\frac{d}{dz} m_{C(w,\de_z), C(1,\de_{\bar{z}^{-1}})}(f)\Big]_{z=1}\Big),$$
$$\tilde{R}_{\de}(f)=2i\pi d(\de)\Big(2m_{C(1,\de), C(1,\de)}(f)+\Big[\frac{d}{dz} (z-1)\Big(m_{P_{\de_z}, C(1,\de_{\bar{z}^{-1}})}(f)+m_{C(w,\de_z), P_{\bar{z}^{-1}}}(f)\Big)\Big]_{z=1}\Big),$$
$$P_\tau(S)=\int_H \textrm{tr} (\tau(h)S) dh,\quad S\in \textrm{End}_{fin.rk}(V_\tau),$$
$$P_{\de_z}(S)=\int_H^* E^0(P,\de_z, S_z)(h) dh,\quad S\in i_{P\cap K}^K \C\otimes  i_{\bar{P}\cap K}^K \check{\C}$$
and 
$$C(s,\de_z)(S):=(1+q^{-1})\int_{K_H\times K_H}\textrm{tr}\big(\big[C^0_{P,P}(s,\de_z)S_z\big](k_1,k_2)\big)dk_1 dk_2,\quad s\in W^G$$

\end{theo}
\section{A local relative trace formula for $PGL(2)$}
We precise the geometric expansion of the truncated kernel  obtained in (\cite{DHSo} Theorem 2.3) for $\sH:=PGL(2)$.  This geometric expansion depends on orbital integrals of $f_1$ and $ f_2$, and on a weight function  $v_L$ where $L=H$ or $M$. To recall the definition of this objects, we need to introduce some notation.\me

\no If $\sJ$ is an algebraic group defined over $\F$, we denote by $J $  its group of points over $\F$ and we identify $\sJ$ with the group of points of $\sJ$ over an algebraic closure of $\F$.  Let $\sJ_H$ be an algebraic subgroup of $\sH$ defined over $\F$.  We denote by  $\sJ:=\mbox{Res}_{\E/\F}(\sJ_H\times_\F\E)$   the restriction of scalars of $\sJ_H$ from $\E$ to $\F$. Then, the group $J:=\sJ(\F)$ is isomorphic to $\sJ_H(\E)$. \me

\no The nontrivial element  of the Galois group of $\E/\F$ 	induces  an involution $\si$ of $\sG$ defined over $\F$. \me

\no We denote by  $\scP$  the connected component of $1$ in the set of $x$ in $\sG$ such that $\si(x)=x^{-1}$. A torus $\sA$ of $\sG$ is called a $\si$-torus if $\sA$ is a torus defined over $\F$ contained in $\scP$. Let $\sS_H$ be a maximal torus of $\sH$.  We denote by   $\sS_\sigma$ the connected component of $\sS\cap \scP$. Then $\sS_\si$ is a maximal $\si$-torus defined over $\F$ and the map $S_H\mapsto S_\si$ is a bijective correspondence  between  $H$-conjugacy classes of   maximal tori of $H$ and $H$-conjugacy classes  of maximal $\si$-tori of $G$. (cf. \cite{DHSo} 1.2).

Each maximal torus of $H$    is either anisotropic   or $H$-conjugate to $M$. We fix $\cT_H$ a set of representatives for the $H$-conjugacy classes of maximal anisotropic torus in $H$.

By (\cite{DHSo} (1.28)), for  each maximal torus $S_H$ of $H$, we can fix a finite set of representatives  $\kappa_S=\{x_m\}$ of the  $(H,S_\si)$-double cosets in $\sH\sS_\si\cap G$ such that each element $x_m$  may be written   $x_m=h_m a^{-1}_m$ where  $h_m\in\sH$  centralizes   the split component $A_{S}$ of $S_H$ and $a_m\in \sS_\si$.\\

The orbital integral of a compactly supported smooth function is defined on the set $G^{\si-reg}$ of $\si$-regular points of $G$, that is the set of point $x$ in $G$ such that $\sH x\sH$ is Zariski closed and of maximal dimension. The set  $G^{\si-reg}$  can be described in terms of maximal $\si$-tori as follows. 
If $\sS_H$ is a maximal torus of $\sH$, we denote by $ \underline{\g s}$ the  Lie algebra of $ \sS $ and we set $\g s := \underline{\g s} (\F)$. We set
$$\De_\si(g)=\textrm{det} (1-\textrm{Ad} (g^{-1}\si(g))_{\g g/\g s }),\quad g\in G .$$ 
By (\cite{DHSo} (1.30)), if  $x\in G^{\si-reg}$  then there exists a maximal torus $S_H$ of $H$ such that $\De_\si(x)\neq 0$. Morever, there are  two elements $x_m\in\kappa_S$ and  $\ga\in S_\si$ such that  $x=x_m\ga$.\\
We define the orbital integral $\cM(f)$  of a function   $f\in C_c^\infty(G)$ on $G^{\si-reg}$ as follows. Let $S_H$ be a maximal torus of $H$. For $x_m\in\kappa_S$ and $\ga\in S_\si$ with $\De_\si(x_m\ga)\neq 0$, we set
\ber $$\cM(f)(x_m\ga):=|\De_\si(x_m\ga)|_\F^{1/4}\int_{ diag(A_S)\bb (H\times H)} f(h^{-1} x_m\ga l) d\overline{(h,l)}$$
where $diag(A_S)$ is the diagonal of $A_S\times A_S$.\eer

We now give an explicit expression of the  truncated function $v_L(\cdot,n)$ defined in (\cite{DHSo} (2.12)), where $n$ is a positive integer and $L$ is   equal to  $H$ or $M$. Let  $n$ be a positive integer. It follows immediately from the definition  (\cite{DHSo} (2.12)) that we have 
\beq  v_H(x_1,y_1,x_2,y_2,n)=1,\quad x_1,y_1,x_2,y_2\in H.\eeq\\
We will describe $v_M$ using (\cite{DHSo} (2.63)). Since $H=P_HK_H$, each  $x\in H$ can be written $x=m_{P_H}(x)n_{P_H}(x)k_{P_H}(x)$ with $m_{P_H}(x)\in M_H, n_{P_H}(x)\in N_H$ and $k_{P_H}(x)\in K_H$. We take similar notation if we consider  $\bar{P}$ instead of $P$. For $Q=P$ or $\bar{P}$, we set 
$$h_{Q_H}(x):=h_{M_H}(m_{Q_H}(x)).$$
With our definition of $h_{M_H}$ (\ref{hM}), the map $M_H\to \R$  given in  (\cite{DHSo} (1.2)) coincides with   $-(\log q) h_{M_H}$ .

For $ x_1,y_1,x_2$ and $y_2$ in $H$, we set 
$$z_P(x_1,y_1,x_2,y_2):=\inf \big( h_{\bar{P}_H}(x_1)-h_{P_H}(y_1),h_{\bar{P}_H}(x_2)- h_{P_H}(y_2)\big),$$
and $$z_{\bar{P}}(x_1,y_1,x_2,y_2):=-\inf\big(h_{\bar{P}_H}(y_1)- h_{P_H}(x_1), h_{\bar{P}_H}(y_2)-h_{P_H}(x_2)\big).$$\\
We omit $ x_1,y_1,x_2$ and $y_2$ in this notation if there is no confusion. Hence  the elements $Z_P^0$ and $Z_{\bar{P}}^0$ of (\cite{DHSo} (2.55)) coincide with  $(\log q)z_P$ and $(\log q)z_{\bar{P}}$ respectively. 
Therefore, the relation  (\cite{DHSo} (2.63)) gives 
$$v_M(x_1,y_1,x_2,y_2,n)=\lim_{\la\to 0} \Big(\frac{q^{\la (n+z_{P})}}{1-q^{-2\la}}(1+q^{-\la})+\frac{q^{\la (-n+z_{\bar{P}})}}{1-q^{2\la}}(1+q^{\la})\Big)$$
$$=\lim_{\la\to 0} \Big(\frac{q^{\la (n+z_{P})}}{1-q^{-\la}} +\frac{q^{-\la (n-z_{\bar{P}})}}{1-q^{ \la}} \Big)=\lim_{\la\to 0}\frac{q^{\la (n+z_{P})}-q^{-\la (n-z_{\bar{P}}+1)}}{1-q^{-\la}} $$
$$=2n+1+ z_P-z_{\bar{P}}.$$

We set $$v^0_M(x_1,y_1,x_2,y_2):= z_P-z_{\bar{P}}$$
$$=\inf \big( h_{\bar{P}_H}(x_1)-h_{P_H}(y_1),h_{\bar{P}_H}(x_2)- h_{P_H}(y_2)\big)+\inf\big(h_{\bar{P}_H}(y_1)- h_{P_H}(x_1), h_{\bar{P}_H}(y_2)-h_{P_H}(x_2)\big).$$

\no Therefore,   (\cite{DHSo} Theorem 2.3) gives:\\
\ber\label{geom}  As $n$ approaches to $+\infty$, the truncated kernel $K^n(f)$ is asymptotic to 
$$2n  \sum_{x_m\in \kappa_M} c^0_{M,x_m}  \int_{M_\si}\cM(f_1)(x_m\ga)\cM(f_2)(x_m\ga) d\ga $$
$$+\sum_{S_H\in\cT_H\cup\{M_H\}}\sum_{x_m\in\kappa_S} c^0_{S,x_m}\int_{S_\si} \cM(f_1)(x_m\ga)\cM(f_2)(x_m\ga) d\ga+  \sum_{x_m\in \kappa_M} c^0_{M,x_m} \int_{M_\si} \cW\cM(f)(x_m\ga)  d\ga,$$\eer
where the constants $ c^0_{M,x_m} $ are defined in (\cite{RR} Theorem 3.4) and $\cW\cM(f)$ is the weighted integral orbital given by $$\De_\si(x_m\ga)^{-1/2}\cW\cM(f)(x_m\ga) $$
$$=\int_{diag(M_H)\bb H\times H}\int_{diag(M_H)\bb H\times H} f_1(x_1^{-1}x_m\ga x_2)  f_2(y_1^{-1} x_m\ga y_2) v^0_M(x_1,y_1,x_2,y_2 )d\overline{(x_1,x_2)} d\overline{(y_1,y_2)}.$$
Therefore, comparing    asymptotic expansions of $K^n(f)$ in  Theorem \ref{SS} and (\ref{geom}), we obtain:
\begin{theo}\label{LRTF} For $f_1$ and $f_2$ in $C_c^\infty(G)$ then we have:
\begin{enumerate}\item 
$$2 \sum_{x_m\in \kappa_M} c^0_{M,x_m}  \int_{M_\si}\cM(f_1)(x_m\ga)\cM(f_2)(x_m\ga) d\ga =\sum_{\de\in \widehat{M}_2,  \de_{|\F^\times}= 1}d(\de)m_{C(1,\de),C(1,\de)}(f).$$
\item (Local relative trace formula). The expression
$$\sum_{S_H\in\cT_H\cup\{M_H\}}\sum_{x_m\in\kappa_S} c^0_{S,x_m}\int_{S_\si} \cM(f_1)(x_m\ga)\cM(f_2)(x_m\ga) d\ga+  \sum_{x_m\in \kappa_M} c^0_{M,x_m} \int_{M_\si} \cW\cM(f)(x_m\ga)  d\ga$$
equals 
$$\sum_{\tau\in \cE_2(G)} d(\tau) m_{P_\tau,P_\tau}(f)+ \frac{1}{4i\pi} \sum_{ \de\in \widehat{M}_2, \de_{|\F^\times}\neq 1, \de_{|\E^1}=1}d(\de) \int_\cO  m_{P_{\de_z},P_{\de_z}}(f)\frac{dz}{z}$$
$$+\frac{1}{4i\pi} \sum_{\de\in \widehat{M}_2, \de_{|\F^\times}= 1} R_{\de}(f)+d(\de)\int_{\cO^-} \frac{ m_{C(1,\de_z), C(1,\de_{\bar{z}^{-1}})}(f)+ m_{C(w,\de_z), C(w,\de_{\bar{z}^{-1}})}(f)}{(1-z)(1-z^{-1})}\frac{dz}{z}$$
$$+\frac{1}{4i\pi} \sum_{\de\in \widehat{M}_2, \de_{|\F^\times}=\de_{|E^1}=1}R_\de(f)+\tilde{R}_{\de}(f)+d(\de)\int_{\cO^-}m_{P_{\de_z},P_{\de_{\bar{z}^{-1}}}}(f)\frac{dz}{z}.$$
where
$$R_{\de}(f):=2i\pi d(\de)\Big(m_{C(1,\de),C(1,\de)}(f)- \Big[\frac{d}{dz} m_{C(w,\de_z), C(1,\de_{\bar{z}^{-1}})}(f)\Big]_{z=1}\Big),$$
$$\tilde{R}_{\de}(f)=2i\pi d(\de)\Big(2m_{C(1,\de), C(1,\de)}(f)+\Big[\frac{d}{dz} (z-1)\Big(m_{P_{\de_z}, C(1,\de_{\bar{z}^{-1}})}(f)+m_{C(w,\de_z), P_{\bar{z}^{-1}}}(f)\Big)\Big]_{z=1}\Big),$$
$$P_\tau(S)=\int_H \textrm{tr} (\tau(h)S) dh,\quad S\in \textrm{End}(V_\tau),$$
$$P_{\de_z}(S)=\int_H^* E^0(P,\de_z, S_z)(h) dh,\quad S\in i_{P\cap K}^K \C\otimes  i_{\bar{P}\cap K}^K \check{\C}$$
and 
$$C(s,\de_z)(S):=(1+q^{-1})\int_{K_H\times K_H}\textrm{tr}\big(\big[C^0_{P,P}(s,\de_z)S_z\big](k_1,k_2)\big)dk_1 dk_2,\quad s\in W^G.$$
\end{enumerate}
\end{theo}
As application of this Theorem, we will invert  orbital integrals on the anisotropic $\si$-torus $M_\si$ of $G$. 

Let $\de\in \widehat{M_2}$. As the operator $C^0_{P,P}(1,\de)$ is the identity operator of $i_P^G\C_{\de_z}\otimes i_{\bar{P}}^G\check{\C_{\de_z}}$, one has 
$$C(1,\de)(v\otimes\check{w})=(1+q^{-1})\int_{K_H\times K_H} v(k_1)\check{w}(k_2) dk_1 dk_2, \quad v\otimes\check{w}\in i_{P\cap K}^K\C\otimes i_{\bar{P}\cap K}^K\check{\C}.$$
Hence, we have $C(1,\de)=(1+q^{-1})\xi_{\de}\otimes \xi_{\check{\de}}$ where $\xi_\de$ and $\xi_{\check{\de}}$ are the $H$-invariant linear forms on 
$ i_{P\cap K}^K\C$ and $i_{\bar{P}\cap K}^K\check{\C}$ respectively given by the integration over $K_H$. Therefore, one deduces that 
$$m_{C(1,\de), C(1,\de)}(f_1\otimes f_2)= m_{\xi_{\de},\xi_{\de}}(f_1)m_{\xi_{\check{\de}},\xi_{\check{\de}}}(f_2).$$
Moreover, by (\cite{AGS} Corollary 5.6.3), the distribution $f\mapsto m_{\xi_{\check{\de}},\xi_{\check{\de}}}(f)$ is smooth in a neighborhood of any $\si$-regular point of $G$. 

\begin{cor} Let $f\in C_c^\infty(G)$. Let $x_m\in\kappa_M$ and $\ga\in M_\si$ such that $x_m\ga$ is $\si$-regular. Then we have
$$c^0_{M,x_m}\cM(f)(x_m\ga)= \sum_{\de\in\widehat{M_2}, \de_{|\F^\times}=1} d(\de) m_{\xi_{\de},\xi_{\de}}(f)m_{\xi_{\check{\de}},\xi_{\check{\de}}}(x_m\ga).$$
\end{cor}
\dem Let $(J_n)_n$ (resp., $(K_n)_n$) be a sequence of compact open sugroups whose intersection is equal to the neutral element of $G$. Then the characteristic function $g_n$ of $J_nx_m\ga K_n$ approaches the Dirac measure at $x_m\ga$. Therefore,  taking  $f_1:=f$ and $f_2:= g_n$ in Theorem \ref{LRTF} {\it1}., we obtain the result.\qed

\no{\bf  Remark.} Let $(\tau, V_\tau)$ be a supercuspidal representation of $G$ and $f$ be a matrix coefficient of $\tau$. Then we deduce from the corollary that the orbital integral of $f$ on $\si$-regular points of $M_\si$ is equal to $0$. 

Moreover,  by (\cite{Fli}, Proposition 11) we have $dim\; V_\tau'^H=1$. Let $\xi$ be a nonzero $H$-invariant linear form on $V_\tau$.  Let $S_H$ be an anisotropic torus of $H$ and $x_m\in\kappa_S$. Then, applying our local relative trace formula to $f_1:=f$ and $f_2$ approaching the Dirac measure at a $\si$-regular point $x_m\ga$ with $\ga\in S_\si$, we obtain
$$\cM(f)(x_m\ga)=c m_{\xi,\xi}(f)m_{\xi,\xi}(x_m\ga),$$
where $c$ is some nonzero constant. 

J. Hakim obtained these results by other methods (\cite{Ha} Proposition 8.1 and Lemma 8.1).


\noindent P.~Delorme, Aix-Marseille Université, CNRS, Centrale Marseille, I2M, UMR 7373, 13453
Marseille, France.\\
{\it E-mail address}:  patrick.delorme@univ-amu.fr\me

\noindent P.~Harinck, CMLS, \'Ecole polytechnique, CNRS-UMR 7640, Université Paris-Saclay, Route de Saclay,  91128 Palaiseau Cedex, France.\\  {\it E-mail address}: pascale.harinck@ polytechnique.edu\me

\end{document}